\documentclass[a4paper]{jpconf}
\usepackage{color} 
\usepackage{xcolor} 
\usepackage{amssymb}
\usepackage{amsmath}
\usepackage{bbold} 
\usepackage{dsfont} 
\usepackage{relsize} 
\usepackage{hyperref}

\def\cmath{\color{black}}
\def\ctxt{\color{black}}
\ifx\alert\undefined\newcommand{\alert}[1]{{\color{red}#1}}\fi
\def\ass{\leftarrow} 
\def\bmat{\begin{pmatrix}}
\def\emat{\end{pmatrix}}
\def\C{\mathbb{C}} 
\def\Q{\mathbb{Q}} 
\def\R{\mathbb{R}} 
\def\Z{\mathbb{Z}}
\newcommand{\baseform}[1]{\mathcal{A}_{#1}} 
\newcommand{\bornform}[1]{\mathcal{B}_{#1}} 
\newcommand{\cabs}[1]{\left|#1\right|} 
\newcommand{\CyclG}[1]{\mathsf{C}_{#1}} 

\newcommand{\farg}[1]{\!\left(#1\right)} 
\newcommand{\GF}[1]{\mathrm{GF}\farg{#1}} 
\DeclareMathOperator{\zeromat}{\mathbb{0}}
\DeclareMathOperator{\idmat}{\mathds{1}}
\def\imu{\mathrm{\mathbf{i}}} 
\newcommand{\IrrRep}[1]{\mathbf{#1}} 
\def\Hspace{\mathcal{H}} 
\newcommand{\Math}[1]{$\cmath{}#1$} 
\newcommand{\MathEq}[1]{\begin{equation*}\cmath{#1}\end{equation*}}
\newcommand{\MathEqLab}[2]{\begin{equation}\cmath{#1}\label{#2}\end{equation}}
\newcommand{\Mone}[1]{\bmat#1\emat} 
\def\NF{\mathcal{F}} 
\def\regrep{\mathrm{P}} 
\def\baseformN{\mathrm{R}} 
\def\repq{\mathsf{U}} 
\DeclareMathOperator{\Ord}{Ord} 
\newcommand{\ordset}[1]{\left[#1\right]} 
\newcommand{\PermRep}[1]{\mathbf{\underline{#1}}} 
\newcommand{\set}[1]{\left\{#1\right\}} 
\newcommand{\transpose}[1]{{#1}^{\mathrm{T}}}
\newcommand{\vect}[1]{\left(#1\right)} 
\def\wG{\mathsf{G}} 
\def\wS{\Omega} 
\def\wSN{\mathsf{N}} 
\begin{document}
\title{A new algorithm for irreducible decomposition of representations of finite groups}
\author{Vladimir V Kornyak}
\address{Laboratory of Information Technologies,
           Joint Institute for Nuclear Research\\
           141980 Dubna, Russia}
\ead{vkornyak@gmail.com}
\begin{abstract}%
An algorithm for irreducible decomposition of representations of finite groups over fields of characteristic zero is described.
The algorithm uses the fact that the 
decomposition induces a partition of the invariant inner product into a complete set of mutually orthogonal projectors.
By expressing the projectors through the basis elements of the centralizer ring of the representation, the problem is reduced to solving systems of quadratic equations.
The current implementation of the algorithm is able to split representations of dimensions up to hundreds of thousands.
Examples of calculations are given.
\end{abstract}
~\\[-50pt]
\section{Introduction}\label{intro} 
The decomposition of linear representations of groups into irreducible subrepresentations is one of the central problems of group theory and its applications in physics.
Currently, the most effective algorithm for solving this problem is a Las Vegas type probabilistic algorithm, called \textbf{\textit{MeatAxe}} \cite{Holt}.
This algorithm is based on the calculation of the characteristic polynomial of a randomly generated matrix of the representation. 
In case of success, factoring this polynomial
and processing the 
factors 
allow either to construct a decomposition of the representation, or to prove its irreducibility.
The \emph{\textbf{MeatAxe}} algorithm played an important role in solving the problem of classifying finite simple groups, where it was applied to representations of groups in linear spaces over small finite fields, such as  \Math{\GF{2}}.
However, \emph{\textbf{MeatAxe}} is inefficient 
in characteristic zero due to the rapid growth of numerical coefficients 
of characteristic polynomials with the matrix dimension, and due to the fact that in characteristic zero a random matrix with high probability has an irreducible characteristic polynomial. 
\par
The quantum formalism is based on Hilbert spaces over fields of characteristic zero.
Traditionally, non-constructive fields \Math{\C} or \Math{\R} are used.
Our goal was to develop an algorithm suitable for the study of quantum-mechanical models based on unitary representations of finite groups over constructive fields of characteristic zero \cite{Kornyak1,Kornyak2}.
The computer implementation of our algorithm, let's call it \textbf{\textit{IrreducibleProjectors}}, splits representations of dimensions up to hundreds of thousands, 
which is not less than the dimensions achievable for \textbf{\textit{MeatAxe}} in the computationally easier context of finite fields.
On the other hand, unlike \textbf{\textit{MeatAxe}}, \textbf{\textit{IrreducibleProjectors}} is of little use in finite-field problems, since it uses the notion of scalar product,
which is problematic for spaces over finite fields.
In fact, \textbf{\textit{IrreducibleProjectors}} and \textbf{\textit{MeatAxe}} have different application areas.
\par
The \textbf{\textit{IrreducibleProjectors}} algorithm requires knowledge of the centralizer ring of the group representation under consideration.
In the general case, the computation of the centralizer ring reduces to a simple problem of linear algebra, namely, to solving a system of matrix equations of the form \Math{AX=XA}.
We will consider here only permutation representations, since (a) any linear representation of a finite group is a subrepresentation of some permutation representation and (b) permutation representations underlie the above mentioned constructive quantum mechanical models.
In the case of permutation representations, the computation of the centralizer ring is particularly simple: it reduces to constructing the orbits of the group on the Cartesian square of the set on which the group acts by permutations.
~\\[-25pt]
\section{Basic concepts and notation}\label{basics} 
Let \Math{\wG} (or, in more detail, \Math{\wG\farg{\wS}}) be a \emph{transitive} permutation group on the set \Math{\wS\cong\set{1,\ldots,\wSN}}.
The action of \Math{g\in\wG} on \Math{i\in\wS} will be denoted by \Math{i^g}.
A \emph{permutation representation} \Math{\regrep} is a representation of \Math{\wG} by matrices of the form \Math{\regrep\farg{g}_{ij}=\delta_{i^{\mathlarger{g}}j}}.
Since \Math{\regrep\farg{g}} is a \Math{\vect{0,1}}-matrix, the permutation representation can be implemented in vector space over any field \Math{\NF}.
We will consider an \Math{\wSN}-dimensional Hilbert space \Math{\Hspace_\wSN} over the field of scalars \Math{\NF}, which is some {constructive} \emph{splitting field} for the group \Math{\wG}.
As \Math{\NF}, one can take a suitable subfield of the \Math{m}th \emph{cyclotomic field}%
, where \Math{m} is the \emph{exponent}%
~of the group \Math{\wG}.
Such a field \Math{\NF}, being an abelian extension of the field of rational numbers \Math{\Q}, is a \emph{constructive dense subfield} of the real \Math{\R} or complex \Math{\C} field.
From the point of view of physics, \Math{\NF} is indistinguishable from \Math{\R} or \Math{\C} and can be freely used in the formalism of quantum mechanics.
\par
An orbit of \Math{\wG} on the Cartesian square \Math{\wS\times\wS} is called an \emph{orbital} \cite{Cameron}.
The number \Math{\baseformN} of orbitals is called the \emph{rank} of the permutation group \Math{\wG\farg{\wS}}.
If the set of orbitals contains some orbital \Math{\Delta}, then it necessarily contains the transposed orbital \Math{\transpose{\Delta}}.
The set of orbitals of a transitive group contains a single \emph{diagonal} orbital \Math{\Delta_1=\set{\vect{i,i}\mid i\in\wS}}, which we will always fix as the first element in the list of orbitals \Math{\set{\Delta_1,\ldots,\Delta_\baseformN}}.
For a transitive group, there is a natural one-to-one correspondence between the orbitals and the orbits of the stabilizer of a point \Math{i\in\wS}, i.e., the subgroup \Math{\wG_i\leq\wG} such that \Math{g\in\wG_i\Rightarrow{}i^g=i}.
An orbit of the stabilizer is called \emph{a suborbit}.
The correspondence between the orbital \Math{\Delta} and the suborbit \Math{\Sigma_i} has the form \Math{\Delta\longleftrightarrow\Sigma_i=\set{j\in\wS\mid\vect{i,j}\in\Delta}.}
The sizes of orbitals and suborbits are related by the equality \Math{\cabs{\Delta}=\wSN\cabs{\Sigma_i}}.
\par
The invariance condition for a bilinear form \Math{A} in the space \Math{\Hspace_\wSN} is expressed by the equations \Math{A=\regrep\farg{g}A\regrep\farg{g^{-1}},~g\in\wG.}
In terms of the matrix entries, these equations have the form \Math{\vect{A}_{\displaystyle{}ij}=\vect{A}_{\displaystyle{}i^gj^g}}.
This implies that the basis of all invariant bilinear forms is in one-to-one correspondence with the set of orbitals.
Namely, to the orbital \Math{\Delta_r\in\set{\Delta_1,\ldots,\Delta_\baseformN}} corresponds the \emph{basis matrix} \Math{\baseform{r}} of size \Math{\wSN\times\wSN} with entries  
\Math{
\vect{\baseform{r}}_{\displaystyle{}ij} =
\begin{cases}
1, &\text{if~}\vect{i,j}\in\Delta_r\,,\\
0, &\text{otherwise}\,.
\end{cases}
}
\par
To implement the algorithms and arrange the output of the results of calculations, it is necessary to introduce some ordering of the basis matrices: 
~\\[-10pt]
\MathEqLab{\baseform{1}\prec\baseform{2}\prec\ldots\prec\baseform{\baseformN}.}{ordbasis}
~\\[-20pt]
We use the following conventions:
\begin{enumerate}
	\item
	\label{inbeg} 
	\Math{\baseform{r}\prec\baseform{s}}, if \Math{\cabs{\Delta_r}<\cabs{\Delta_s}} (or, equivalently,  \Math{\cabs{(\Sigma_i)_r}<\cabs{(\Sigma_i)_s}} --- comparing suborbit lengths),
	\item \Math{\baseform{r}\prec\baseform{s}}, if \Math{\baseform{r}=\transpose{\mathcal{A}}_r\wedge\baseform{s}\neq\transpose{\mathcal{A}}_s} (symmetric matrices precede asymmetric),
	\item
		\label{1stcolumn} \Math{\baseform{r}\prec\baseform{s}}, if \Math{{I}_{\baseform{r}}<{I}_{\baseform{s}}}, where \Math{I_X=\min\vect{i\mid(X)_{i1}=1}} (comparing the positions of the first nonzero element in the first columns of matrices),
	\item
	\label{enbeg} if \Math{\baseform{r}\neq\transpose{\mathcal{A}}_r}, then \Math{\baseform{r+1}=\transpose{\mathcal{A}}_r} (paired matrices are always placed adjacently).
\end{enumerate}
Applying rules \ref{inbeg} --- \ref{enbeg} in the specified order uniquely defines the sequence \eqref{ordbasis}.
According to these rules, the diagonal orbital matrix is the first element of the list \eqref{ordbasis}:  \Math{\baseform{1}=\idmat_{\wSN}}.
\par
The set of invariant bilinear forms has the structure of a ring, which is called the \emph{centralizer ring} (or \emph{centralizer algebra}).
The multiplication table for basic elements \eqref{ordbasis} has the form
~\\[-10pt]
\MathEqLab{\hspace*{120pt}\baseform{p}\baseform{q}=\sum_{r=1}^{\baseformN}C_{pq}^r\baseform{r},}{multtab}
~\\[-10pt]
where the coefficients \Math{C_{pq}^r} are natural numbers lying within \Math{0\leq{}C_{pq}^r<\wSN}.
The representation \Math{\regrep} is \emph{multiplicity-free} if and only if the centralizer ring is \emph{commutative}.
~\\[-25pt]
\section{Algorithm description}\label{algo} 
Let \Math{T} be a unitary (we can always provide unitarity) transformation matrix splitting the representation \Math{\regrep} in the Hilbert space \Math{\Hspace_\wSN} into \Math{M} irreducible components:
~\\[-10pt]
\MathEq{T^{-1}\regrep\farg{g}T=1\oplus\repq_{\!d_2}\farg{g}\oplus\cdots\oplus\repq_{\!d_m}\farg{g}\oplus\cdots\oplus\repq_{\!d_M}\farg{g},}
~\\[-15pt]
where \Math{\repq_{\!d_m}} is  a \Math{d_m}-dimensional irreducible component.
\par
The \emph{standard scalar product} in the Hilbert space is represented by the matrix  \Math{\idmat_\wSN} in any orthonormal basis.
In the splitting basis, we have the following decomposition
~\\[-10pt]
\MathEqLab{\idmat_\wSN=\idmat_{d_1=1}\oplus\cdots\oplus\idmat_{d_m}\oplus\cdots\oplus\idmat_{d_M}.}{idsplit}
~\\[-15pt]
Here \Math{\idmat_{d_1=1}\equiv\Mone{1}} is the scalar product in the one-dimensional \emph{trivial subrepresentation} that is always present in any permutation representation. 
The \emph{preimage} of decomposition \eqref{idsplit} in the original permutation basis has the form
~\\[-10pt]
\MathEqLab{\idmat_\wSN=\bornform{1}+\cdots+\bornform{m}+\cdots+\bornform{M},}{completeness}
~\\[-15pt]
where \Math{\bornform{m}} is defined by the relation
~\\[-10pt]
\MathEqLab{T^{-1}\bornform{m}T=\zeromat_{1+d_2+\cdots+d_{m-1}}\oplus\idmat_{d_m}\oplus\zeromat_{d_{m+1}+\cdots+d_M}\equiv\mathcal{D}_m.}{brel}
~\\[-15pt]
It can be seen from this relation that the matrices  \Math{\bornform{m}} are \emph{idempotent}
~\\[-10pt]
\MathEqLab{\bornform{m}^2=\bornform{m}}{idempotency}
~\\[-20pt]
and \emph{mutually orthogonal}
~\\[-15pt]
\MathEqLab{\bornform{m}\bornform{m'}=\zeromat_\wSN~\text{\ctxt{}if}~m\neq{m'}.}{orthogonality}
~\\[-15pt]
Relations \eqref{idempotency} and \eqref{orthogonality} together with the \emph{completeness condition} \eqref{completeness} mean that the set \Math{\bornform{1},\ldots,\bornform{M}} is a \emph{complete system of mutually orthogonal projectors} in the Hilbert space \Math{\Hspace_\wSN}.
\par
The set of \emph{irreducible invariant} projectors \Math{\bornform{1},\ldots,\bornform{M}} contains complete information about the decomposition of the representation \Math{\regrep} into irreducible components.
In particular, the transformation matrix \Math{T} can be computed by solving the system of linear equations
~\\[-10pt]
\MathEq{\bornform{1}T-T\mathcal{D}_1=\cdots=\bornform{M}T-T\mathcal{D}_M=\zeromat_\wSN.}
Any invariant projector is a solution of the equation
~\\[-10pt]
\MathEqLab{X^2-X=\zeromat_\wSN,}{ideq}
~\\[-15pt]
where \Math{X=x_1\baseform{1}+\cdots+x_\baseformN\baseform{\baseformN}} is a generic invariant bilinear form written in basis \eqref{ordbasis}.
Using multiplication table \eqref{multtab} and decomposing \eqref{ideq} into components in basis (\eqref{ordbasis}, we obtain the system of \Math{\baseformN} quadratic equations for \Math{\baseformN} unknowns \Math{x_1,\ldots,x_\baseformN}
~\\[-10pt]
\MathEqLab{E\farg{x_1,\ldots,x_\baseformN}=0\sim\set{E_1\farg{x_1,\ldots,x_\baseformN}=0,\ldots,E_\baseformN\farg{x_1,\ldots,x_\baseformN}=0}.}{ideqpoly}
~\\[-15pt]
We will call the left hand sides of these equations \emph{idempotency polynomials}.
An irreducible invariant projector \Math{\bornform{m}} in basis has the form
~\\[-10pt]
\MathEqLab{\bornform{m}=b_{m,1}\baseform{1}+b_{m,2}\baseform{2}+\cdots+b_{m,\baseformN}\baseform{\baseformN},}{bmina}
~\\[-15pt]
where the vector \Math{B_m=\ordset{b_{m,1},\ldots,b_{m,\baseformN}}} is a solution of the system of equations \eqref{ideqpoly}.
Due to the invariance of the trace of a matrix under the similarity transformation, relation \eqref{brel} implies the equality \Math{\tr\bornform{m}=d_{m}.}
Combining this equality with the fact that in \eqref{bmina} only \Math{\baseform{1}} has nonzero diagonal elements and \Math{\tr\baseform{1}=\wSN},
we can fix the first coefficient in decomposition \eqref{bmina}: 
~\\[-10pt]
\MathEq{b_{m,1}=d_m/\wSN.}
~\\[-15pt]
Thus, the possible values of \Math{x_1} that provide solutions of the polynomial system \eqref{ideqpoly} are fractions of the form \Math{d/\wSN},
where natural numbers \Math{d\in\ordset{1,\ldots,\wSN-1}} are either irreducible dimensions \Math{d_m} or sums of such dimensions.
Orthogonality condition \eqref{orthogonality} allows us to exclude from consideration dimensions that are not irreducible.
For generic \Math{B=b_1\baseform{1}+\cdots+b_\baseformN\baseform{\baseformN}} and \Math{X}, the orthogonality condition can be written as 
~\\[-15pt]
\MathEqLab{BX=\zeromat_\wSN.}{genort}
~\\[-15pt]
This matrix equation is a system of linear with respect to variables \Math{x_1,\ldots,x_\baseformN} equations with parameters \Math{b_1,\ldots,b_\baseformN}.   
Using multiplication table \eqref{multtab}, the left hand side of \eqref{genort} can be represented as a system of \Math{\baseformN} bilinear forms
~\\[-10pt]
\MathEqLab{O\farg{b_1,\ldots,b_\baseformN;x_1,\ldots,x_\baseformN}=
\left\{\hspace*{-10pt}\text{
\begin{tabular}{c}
$O_1\farg{b_1,\ldots,b_\baseformN;x_1,\ldots,x_\baseformN},$\\
$\vdots$\\
$O_{\baseformN}\farg{b_1,\ldots,b_\baseformN;x_1,\ldots,x_\baseformN}.$
\end{tabular}
\hspace*{-10pt}}\right\},
}{orthpoly}
~\\[-10pt]
which we will call \emph{orthogonality polynomials}.
\par
The main part of the algorithm is organized as a cycle starting with \Math{d=1} and ending when the sum of the irreducible dimensions reaches the value \Math{\wSN}.
The current \Math{d} is processed as follows:
~\\[-15pt]
\begin{enumerate}
	\item 
Substitute \Math{x_1=d/\wSN} into \eqref{ideqpoly} and solve the system of equations 
~\\[-10pt]
\MathEqLab{E\farg{d/\wSN,x_2,\ldots,x_\baseformN}=0.}{trunkeq}
~\\[-15pt]
At the same time, without significant additional calculations, the Hilbert dimension \Math{h} of the corresponding polynomial ideal is determined.
The solution is always realizable algorithmically, since all the roots of the system belong to abelian extensions of \Math{\Q}.
Modern computer algebra systems, in particular \textbf{\textit{Maple}}, cope well with this task.
	\item 
If system \eqref{trunkeq} is incompatible, then the current value of \Math{d} is not an irreducible dimension and we go to the next value of \Math{d} in the loop.	
	\item \label{h=0}
If the Hilbert dimension \Math{h=0} and system \eqref{trunkeq} has \Math{k} solutions, then we get \Math{k} (different if \Math{k>1}) \Math{d}-dimensional irreducible subrepresentations.	
	\item \label{h>0}
\Math{h>0} indicates a \Math{d}-dimensional irreducible component of the nontrivial multiplicity \Math{k}.
The corresponding component of the centralizer ring has the structure  \Math{A\otimes\idmat_d}, where \Math{A} is an arbitrary matrix of size \Math{k\times{}k}.
The idempotency condition, \Math{\vect{A\otimes\idmat_d}^2=A\otimes\idmat_d,} imposes the constraint on \Math{A}:
\Math{A^2-A=0.}
The complete family of solutions of this equation %
is a manifold of dimension \Math{{h}=\left\lfloor{{k}^2}/{2}\right\rfloor.}
Hence, for the multiplicity, we have: \Math{{k}=\left\lceil\sqrt{2h}\right\rceil.}
\par
Then, using some procedure, \Math{k} arbitrary but mutually orthogonal representatives are selected from the family of equivalent \Math{d}-dimensional projectors.
	\item \label{prjproc}
Each of the k irreducible projectors obtained in items \eqref{h=0} or \eqref{h>0} is processed as follows.
Projector \Math{\bornform{m}} is added to the list of irreducible projectors.
The corresponding invariant subspace is excluded from further consideration by adding the orthogonality polynomials \Math{\bornform{m}X} to the set of polynomials \eqref{ideqpoly}:
\Math{E\farg{x_1,x_2,\ldots,x_\baseformN}\ass{}E\farg{x_1,x_2,\ldots,x_\baseformN}\cup\set{\bornform{m}X}.}
	\item 
After the described in item \eqref{prjproc} processing of all \Math{k} irreducible projectors, the transition to the next \Math{d} is performed.	
\end{enumerate}
\par 
The \textbf{\textit{IrreducibleProjectors}} algorithm is implemented in the form of two procedures, called \textbf{\textit{PreparePolynomialData}} and \textbf{\textit{SplitRepresentation}}.
\begin{enumerate}
	\item 
The \textbf{\textit{PreparePolynomialData}} procedure is implemented in \textbf{\textit{C}}.	
The input is the set of generators of \Math{\wG\farg{\wS}.}
The program computes basis \eqref{ordbasis}, multiplication table \eqref{multtab}, constructs polynomials of idempotency \eqref{ideqpoly} and orthogonality \eqref{orthpoly},
and the code for the procedure \textbf{\textit{SplitRepresentation}}.
This code is task-specific: for non-commutative centralizer ring some additional functions to process multiple subrepresentations are generated.
	\item 
\textbf{\textit{SplitRepresentation}} is a \textbf{\textit{Maple}} code generated by the \textbf{\textit{PreparePolynomialData}}.
This code performs the above-described cycle over dimensions.
The polynomial systems are processed by functions from the \textbf{\textit{Groebner}} package implemented in \textbf{\textit{Maple}}.
\end{enumerate}
Algorithms and related implementation and technical issues are described in more detail in \cite{Kornyak3}.
\section{Examples of calculations}\label{calc} 
The input data are taken from the ``Sporadic groups'' section of the \textsc{Atlas}  \cite{atlas}.
The \textsc{Atlas} contains representations of simple groups and some of their extensions.
Namely, if a group \Math{\wG} has a non-trivial
~\\[-15pt]
\begin{enumerate}
	\item 
second homology group \Math{H_2\farg{\wG,\Z}}, called the \emph{Schur multiplier} and denoted by the symbol \Math{\mathrm{M}\farg{\wG}}, then there are nontrivial central extensions of \Math{\wG} by subgroups of \Math{\mathrm{M}\farg{\wG}}; 
	\item 
\emph{outer automorphisms} \Math{\mathrm{Out}\farg{\wG}}, 
then there are extensions with \Math{\wG} as a normal subgroup.
\end{enumerate}
~\\[-15pt]
\Math{A.B} denotes a generic extension of \Math{B} by \Math{A}.  
A \emph{split} extension is denoted by \Math{A\rtimes{}B}.
Cyclic groups \Math{\CyclG{n}} are represented by their orders \Math{n} in the notation for extensions.
\par
We have tried for completeness to choose examples from all generations of the ``Happy Family'' and from the ``Pariahs'' family.
\par
Irreducible components are denoted by their dimensions in bold (possibly with additional indices to distinguish between non-equivalent subrepresentations of the same dimension).
Permutation representations are denoted by their dimensions in bold with an underscore.
\Math{\bornform{\IrrRep{m}}} denotes  the irreducible projector corresponding to the irreducible subrepresentation \Math{\mathbf{m}}.
\par
The calculations were performed on a PC with a 3.30GHz 
CPU and 16GB RAM.
\vspace*{-10pt}
\subsection{Detailed example}\label{example} 
Here is a compact example of the outputs produced by the programs.
The \textbf{Held group} \Math{\mathbf{He}} has the properties:
\Math{\Ord\farg{H\!e}=4030387200=2^{10}\cdot3^3\cdot5^2\cdot7^3\cdot17,~ \mathrm{M}\farg{H\!e}\cong1,}~ \Math{\mathrm{Out}\farg{H\!e}\cong\CyclG{2}.}
\par
The program \textbf{\textit{PreparePolynomialData}}, applied to the \textbf{8330}-dimensional representation of this group, in addition to the code of the program \textbf{\textit{SplitRepresentation}} and input data for it, produces the following text:
\begin{verbatim}
___Action of He on 8330 points
Rank of He_on_8330: 7
Dimension: 8330
Suborbit lengths: 1, 105, 720, 840, 840', 1344, 4480.
Centralizer ring is commutative
=> permutation representation is multiplicity free
___Total time: 2.93 sec
___Technical information
Orbital matrices space: 57.9 MB
Orbital path space    : 35.6 MB
Total orbital space   : 93.5 MB
Maximum number of polynomial terms: 217
\end{verbatim}
This text contains information about the rank of the representation, the lengths of the suborbits (the pair \verb"840, 840'" refers to the mutually transposed orbitals), the presence or absence of multiple subrepresentations, as well as the time and memory spent to solve the problem.
\par
\textit{\textbf{SplitRepresentation}} produces the following decomposition and invariant projectors
~\\[-10pt]
\MathEq{\PermRep{8330}\cong\IrrRep{1}\oplus\IrrRep{51}\oplus\overline{\IrrRep{51}}\oplus\IrrRep{680}\oplus\IrrRep{1275}\oplus\IrrRep{1920}\oplus\IrrRep{4352}}
~\\[-30pt]
\begin{align*}
\bornform{\IrrRep{1}} =&~ \frac{1}{8330}\left(\baseform{1}+\baseform{2}+\baseform{3}+\baseform{4}+\baseform{5}+\baseform{6}+\baseform{7}\right)\\
\bornform{\IrrRep{51}} =&~ \frac{3}{490}\left(\baseform{1}+\frac{\baseform{2}}{3}-\frac{\baseform{3}}{6}-\frac{1-\imu\sqrt{7}}{12}\baseform{4}-\frac{1+\imu\sqrt{7}}{12}\baseform{5}+\frac{\baseform{6}}{6}\right)\\
\bornform{\IrrRep{680}} =&~ \frac{4}{49}\left(\baseform{1}+\frac{\baseform{2}}{5}+\frac{\baseform{3}}{120}+\frac{\baseform{4}}{20}+\frac{\baseform{5}}{20}-\frac{\baseform{7}}{40}\right)\\
\end{align*}
\begin{align*}
\bornform{\IrrRep{1275}} =&~ \frac{15}{98}\left(\baseform{1}+\frac{\baseform{2}}{15}+\frac{\baseform{3}}{15}-\frac{\baseform{4}}{30}-\frac{\baseform{5}}{30}\right)\\
\bornform{\IrrRep{1920}} =&~ \frac{192}{833}\left(\baseform{1}-\frac{2\baseform{2}}{15}+\frac{\baseform{3}}{120}+\frac{\baseform{4}}{120}+\frac{\baseform{5}}{120}+\frac{5\baseform{6}}{192}-\frac{3\baseform{7}}{320}\right)\\
\bornform{\IrrRep{4352}} =&~ \frac{128}{245}\left(\baseform{1}-\frac{\baseform{3}}{48}-\frac{\baseform{6}}{64}+\frac{\baseform{7}}{128}\right)
\end{align*}
~\\[-15pt]
\texttt{Time: 1.4 sec}\\
Here \Math{\IrrRep{51}} and \Math{\overline{\IrrRep{51}}} are two different complex conjugate representations of dimension 51.
~\\[-20pt]
\subsection{Comparison with the implementation of \textbf{\textit{MeatAxe}} in \textbf{\textit{Magma}}}\label{vsmagma} 
The \textbf{\textit{Magma}}  implementation of the \textbf{\textit{MeatAxe}} algorithm is considered one of the best.
The \textbf{\textit{Magma}} database contains a 3906-dimensional permutation representation of the group \Math{G_2\farg{5}} -- an exceptional group of Lie type.
The decomposition of this representation into irreducible components over the field \Math{\GF{2}} is given in \cite{Magma} to illustrate the possibilities of \textit{\textbf{MeatAxe}}.
\par
The application of our programs to this representation gives the following data:
\par
Rank: \Math{4}. Suborbit lengths: \Math{1, 30, 750, 3125}.
\begin{align*}
\PermRep{3906}\cong&~ \IrrRep{1}\oplus\IrrRep{930}\oplus\IrrRep{1085}\oplus\IrrRep{1890}
\\
    \bornform{\IrrRep{1}} =&~ \frac{1}{3906}\vect{\baseform{1}+\baseform{2}+\baseform{3}+\baseform{4}}
\\
\bornform{\IrrRep{930}} =&~ \frac{5}{21}\vect{\baseform{1}+\frac{3}{10}\baseform{2}+\frac{1}{50}\baseform{3}-\frac{1}{125}\baseform{4}}
\\
\bornform{\IrrRep{1085}} =&~ \frac{5}{18}\vect{\baseform{1}-\frac{1}{5}\baseform{2}+\frac{1}{25}\baseform{3}-\frac{1}{125}\baseform{4}}
\\
\bornform{\IrrRep{1890}} =&~ \frac{15}{31}\vect{\baseform{1}-\frac{1}{30}\baseform{2}-\frac{1}{30}\baseform{3}+\frac{1}{125}\baseform{4}}
\end{align*}
\par
Time \textbf{C}: \Math{0.5} sec. Time \textbf{\textbf{Maple}}: \Math{0.8} sec.\\
We see that in the characteristic zero the representation splits over the field \Math{\Q}.
\par
Splitting this representation over \Math{\Q} using \textbf{\textit{Magma}} fails due to memory exhaustion.
However, it is possible to reproduce the same set of irreducible dimensions as in the case of characteristic zero, if we split the representation over a finite field with a characteristic that does not divide the order of the group.
In our case, we have \Math{\Ord\farg{{G_2\farg{5}}}=5859000000=2^6\cdot3^3\cdot5^6\cdot7\cdot31}.
Therefore, the smallest field that ``mimics'' \Math{\Q} in the above sense is \Math{\GF{11}}.
We present a session of the corresponding computation using \textbf{\textit{Magma}} (execution time is given in seconds).
\begin{verbatim}
> load "g25";
Loading "/opt/magma.21-1/libs/pergps/g25"
The Lie group G( 2, 5 ) represented as a permutation
group of degree 3906.
Order: 5 859 000 000 = 2^6 * 3^3 * 5^6 * 7 * 31.
Group: G
> time Constituents(PermutationModule(G,GF(11)));
[
    GModule of dimension 1 over GF(11),
    GModule of dimension 930 over GF(11),
    GModule of dimension 1085 over GF(11),
    GModule of dimension 1890 over GF(11)
]
Time: 282.060
\end{verbatim}
~\\[-45pt]
\subsection{Some calculations for sporadic groups}\label{sporadic} 
\par
The data below contain information about ranks, suborbit lengths, structures of irreducible decompositions, and calculation times.
For brevity, we omitted explicit expressions for irreducible projectors  \Math{\bornform{\IrrRep{m}}}.
The expression \Math{\ell^m} in the list of suborbit lengths means that there are \Math{m} suborbits of length \Math{\ell}.
Non-equivalent irreducible components of the same dimension differ, either by the symbol of complex conjugation (overbar), or by the Greek indices, or by the indices \Math{\pm}, meaning that there are two components having the structure \Math{A\pm{}B}.
Multiple subrepresentations are underbraced.
The execution times are given separately for \textit{\textbf{PreparePolynomialData}} (Time \textbf{C}) and \textit{\textbf{SplitRepresentation}} (Time \textbf{Maple}).
~\\[-25pt]
\subsubsection{\textbf{Mathieu groups.}}\label{Mathieu} 
The five Mathieu groups \Math{M_{11}}, \Math{M_{12}}, \Math{M_{22}}, \Math{M_{23}} and \Math{M_{24}} are the first sporadic groups that have been discovered.
Each group \Math{M_{n}} is isomorphic to a \emph{multiply transitive} permutation group on \Math{n} elements.
The 5-transitive group \Math{M_{12}} and the 3-transitive group \Math{M_{22}} are the only Mathieu groups that have non-trivial Schur multipliers and outer automorphism groups.
From the point of view of the structure of irreducible decompositions, the most interesting are the \emph{covers} of the \textbf{Mathieu group} \Math{\mathbf{M_{22}}}.\\
Main properties of \Math{M_{22}}:
\Math{\Ord\farg{M_{22}}=443520=2^{7}\cdot3^2\cdot5\cdot7\cdot11,~ \mathrm{M}\farg{M_{22}}\cong\CyclG{12},~ \mathrm{Out}\farg{M_{22}}\cong\CyclG{2}.}
~\\[-20pt]
\begin{enumerate}
	\item 
\textbf{990}-dimensional representation of \Math{\mathbf{3.M_{22}}}\\
Rank: \Math{13}. Suborbit lengths: \Math{1^3, 7^3, 42^3, 168^3, 336}.\\[1pt]
{\small\Math{\PermRep{990}\cong\IrrRep{1}\oplus\IrrRep{21_\alpha}\oplus\IrrRep{21_\beta}\oplus\overline{\IrrRep{21_\beta}}\oplus\IrrRep{55}
\oplus\IrrRep{99_\alpha}\oplus\IrrRep{99_\beta}\oplus\overline{\IrrRep{99_\beta}}
\oplus\IrrRep{105_+}\oplus\overline{\IrrRep{105_+}}\oplus\IrrRep{105_-}\oplus\overline{\IrrRep{105_-}}\oplus\IrrRep{154}}}\\[1pt]
Time \textbf{C}: \Math{1} sec. Time \textbf{\textbf{Maple}}: \Math{28} sec.
	\item 
~\\[-17pt]
\textbf{2016}-dimensional representation of \Math{\mathbf{3.M_{22}}}\\
Rank: \Math{16}. Suborbit lengths: \Math{1^3, 55^3, 66^3, 165^4, 330^3}.\\[1pt]
{\small\Math{\PermRep{2016}\cong\IrrRep{1}\oplus\IrrRep{21_\alpha}\oplus\IrrRep{21_\beta}\oplus\overline{\IrrRep{21_\beta}}\oplus\IrrRep{55}
\oplus\IrrRep{105_+}\oplus\overline{\IrrRep{105_+}}\oplus\IrrRep{105_-}\oplus\overline{\IrrRep{105_-}}}\\
\Math{\hspace{42pt}\oplus\IrrRep{154}\oplus\IrrRep{210_\alpha}\oplus\IrrRep{210_\beta}\oplus\overline{\IrrRep{210_\beta}}\oplus\IrrRep{231_\alpha}\oplus\IrrRep{231_\beta}\oplus\overline{\IrrRep{231_\beta}}}}\\[1pt]
Time \textbf{C}: \Math{2} sec. Time \textbf{\textbf{Maple}}: \Math{1} h \Math{15} min \Math{52} sec.
	\item
~\\[-17pt]
\textbf{1980}-dimensional representation of \Math{\mathbf{6.M_{22}}}\\
Rank: \Math{17}. Suborbit lengths: \Math{1^6, 14^3, 84^3, 336^5}.\\[1pt]
{\small\Math{\PermRep{1980}\cong\IrrRep{1}\oplus\IrrRep{21_\alpha}\oplus\IrrRep{21_\beta}\oplus\overline{\IrrRep{21_\beta}}\oplus\IrrRep{55}\oplus\IrrRep{99_\alpha}\oplus\IrrRep{99_\beta}\oplus\overline{\IrrRep{99_\beta}}\oplus\IrrRep{105_+}\oplus\overline{\IrrRep{105_+}}\oplus\IrrRep{105_-}\oplus\overline{\IrrRep{105_-}}}\\
\Math{\hspace{42pt}\oplus\IrrRep{120}\oplus\IrrRep{154}\oplus\IrrRep{210}\oplus\IrrRep{330}\oplus\overline{\IrrRep{330}}}}\\[1pt]
Time \textbf{C}: \Math{1} sec. Time \textbf{\textbf{Maple}}: \Math{6} h \Math{34} min \Math{14} sec.
\end{enumerate}
~\\[-43pt]
\subsubsection{\textbf{Leech lattice groups.}}\label{Leech} 
~\\[-27pt]
\paragraph{\textbf{Higman-Sims group} \Math{\textbf{HS}}.}
 \Math{\Ord=44352000=2^{9}\cdot3^2\cdot5^3\cdot7\cdot11, ~\mathrm{M}\cong\CyclG{2},~ \mathrm{Out}\cong\CyclG{2}.}
~\\[-20pt]
\begin{enumerate}
	\item 
\textbf{5600}-dimensional representation of \Math{\mathbf{HS}}\\
Rank: \Math{9}. Suborbit lengths: \Math{1, 55, 132, 165, 495, 660, 792, 1320, 1980}.\\[1pt]
{\small\Math{\PermRep{5600}\cong\IrrRep{1}\oplus\IrrRep{22}\oplus\IrrRep{77}\oplus\IrrRep{154}\oplus\IrrRep{175}\oplus\IrrRep{770}
\oplus\IrrRep{825}\oplus\IrrRep{1056}\oplus\IrrRep{2520}}}\\[1pt]
Time \textbf{C}: \Math{2} sec. Time \textbf{\textbf{Maple}}: \Math{2} sec.
	\item 
~\\[-17pt]
\textbf{11200}-dimensional representation of \Math{\mathbf{2.HS}}\\
Rank: \Math{16}. Suborbit lengths: \Math{1^2, 110, 132^2, 165^2, 660^2, 792^2, 990, 1320^2, 1980^2}.\\[1pt]
{\small\Math{\PermRep{11200}\cong\IrrRep{1}\oplus\IrrRep{22}\oplus\IrrRep{56}\oplus\IrrRep{77}\oplus\IrrRep{154}
\oplus\IrrRep{175}  \oplus\IrrRep{176}\oplus\overline{\IrrRep{176}}\oplus\IrrRep{616}\oplus\overline{\IrrRep{616}}
\oplus\IrrRep{770}\oplus\IrrRep{825}\oplus\IrrRep{1056}\oplus\IrrRep{1980}\oplus\overline{\IrrRep{1980}}\oplus\IrrRep{2520}}}\\[1pt]
Time \textbf{C}: \Math{7} sec. Time \textbf{\textbf{Maple}}: \Math{1} h \Math{25} min \Math{47} sec.
	\item 
~\\[-17pt]
\textbf{1100}-dimensional representation of \Math{\mathbf{HS\rtimes{2}}}\\
Rank: \Math{5}. Suborbit lengths: \Math{1, 28, 105, 336, 630}.\\[1pt]
{\small\Math{\PermRep{1100}\cong\IrrRep{1}\oplus\IrrRep{77}\oplus\IrrRep{154}\oplus\IrrRep{175}\oplus\IrrRep{693}}}\\[1pt]
Time \textbf{C}: \Math{<1} sec. Time \textbf{\textbf{Maple}}: \Math{<1} sec.
	\item 
~\\[-17pt]
\textbf{1408}-dimensional representation of \Math{\mathbf{2.HS.2}}\\
Rank: \Math{11}. Suborbit lengths: \Math{1^4, 50^4, 350^2, 504}.\\[1pt]
{\small\Math{\PermRep{1408}\cong\IrrRep{1}\oplus\IrrRep{1'}\oplus\IrrRep{22_+}\oplus\IrrRep{22_-}\oplus\IrrRep{175_+}\oplus\IrrRep{175_-}\oplus\IrrRep{308}
\oplus\underbrace{\IrrRep{352}\oplus\IrrRep{352}}_{}}}\\[-9pt]
Time \textbf{C}: \Math{<1} sec. Time \textbf{\textbf{Maple}}: \Math{3} sec.
\end{enumerate}
~\\[-40pt]
\paragraph{\textbf{Janko group} \Math{\mathbf{J_2}}.}\label{janko}
\Math{\Ord=604800=2^{7}\cdot3^3\cdot5^2\cdot7, \mathrm{M}\cong\CyclG{2}, \mathrm{Out}\cong\CyclG{2}.}\\
\textbf{1800}-dimensional representation of \Math{\mathbf{J_2}}\\
Rank: \Math{18}. Suborbit lengths: \Math{1, 14^2, 21, 28, 42^3, 84^3, 168^6, 336}.\\[1pt]
{\small\Math{\PermRep{1800}\cong\IrrRep{1}\oplus\IrrRep{36}\oplus\underbrace{\IrrRep{63}\oplus\IrrRep{63}}_{}\oplus\underbrace{\IrrRep{126}\oplus\IrrRep{126}}_{}\oplus\IrrRep{160}\oplus\IrrRep{175}\oplus\IrrRep{288}\oplus\underbrace{\IrrRep{336}\oplus\IrrRep{336}}_{}}}\\[-5pt]
Time \textbf{C}: \Math{2} sec. Time \textbf{\textbf{Maple}}: \Math{13} min \Math{29} sec.
~\\[-27pt]
\paragraph{\textbf{Conway group} \Math{\mathbf{Co_1}}.\hspace*{-5pt}}\label{Co_1}
\Math{\Ord=4157776806543360000=2^{21}\cdot3^9\cdot5^4\cdot7^2\cdot11\cdot13\cdot23,~ \mathrm{M}\cong\CyclG{2},~ \mathrm{Out}\cong1.}\\
\textbf{98280}-dimensional representation of \Math{\mathbf{Co_1}}\\
Rank: \Math{4}. Suborbit lengths: \Math{1, 4600, 46575, 47104}.\\
{\small\Math{\PermRep{98280}\cong\IrrRep{1}\oplus\IrrRep{299}\oplus\IrrRep{17250}\oplus\IrrRep{80730}}}\\[1pt]
Time \textbf{C}: \Math{43} min \Math{12} sec. Time \textbf{\textbf{Maple}}: \Math{6} sec.\\
\emph{Remark.} The program \textit{\textbf{PreparePolynomialData}} uses more than 8.8 GB of RAM for this task.
~\\[-27pt]
\paragraph{\textbf{Conway group} \Math{\mathbf{Co_2}}.}\label{Co_2}
\Math{\Ord=42305421312000=2^{18}\cdot3^6\cdot5^3\cdot7\cdot11\cdot23,~ \mathrm{M}\cong1,~ \mathrm{Out}\cong1.}\\
\textbf{4600}-dimensional representation of \Math{\mathbf{Co_2}}\\
Rank: \Math{5}. Suborbit lengths: \Math{1^2, 891^2, 2816}.\\[1pt]
{\small\Math{\PermRep{4600}\cong\IrrRep{1}\oplus\IrrRep{23}\oplus\IrrRep{275}\oplus\IrrRep{2024}\oplus\IrrRep{2277}}}\\[1pt]
Time \textbf{C}: \Math{<1} sec. Time \textbf{\textbf{Maple}}: \Math{<1} sec.
~\\[-27pt]
\paragraph{\textbf{Conway group} \Math{\mathbf{Co_3}}.}\label{Co_3}
\Math{\Ord=495766656000=2^{10}\cdot3^7\cdot5^3\cdot7\cdot11\cdot23,~ \mathrm{M}\cong1,~ \mathrm{Out}\cong1.}\\
\textbf{48600}-dimensional representation of \Math{\mathbf{Co_3}}\\
Rank: \Math{8}. Suborbit lengths: \Math{1, 253, 506, 1771, 7590, 8855, 14168, 15456}.\\[1pt]
{\small\Math{\PermRep{48600}\cong\IrrRep{1}\oplus\IrrRep{23}\oplus\IrrRep{253}\oplus\IrrRep{275}\oplus\IrrRep{2024}\oplus\IrrRep{5544}\oplus\IrrRep{8855}\oplus\IrrRep{31625}}}\\[1pt]
Time \textbf{C}: \Math{2} min \Math{17} sec. Time \textbf{\textbf{Maple}}: \Math{2} sec.
~\\[-27pt]
\paragraph{\textbf{McLaughlin group} \Math{\mathbf{McL}}.}\label{McL}
\Math{\Ord=898128000=2^{7}\cdot3^6\cdot5^3\cdot7\cdot11,~ \mathrm{M}\cong\CyclG{3},~ \mathrm{Out}\cong\CyclG{2}.}
~\\[-19pt]
\begin{enumerate}
	\item 
\textbf{22275}-dimensional representation \textbf{(a)} of \Math{\mathbf{McL}}\\
Rank: \Math{13}. Suborbit lengths: \Math{1, 112, 140, 210, 420, 672, 1680^2, 2240, 3360^3, 5040}.\\[1pt]
{\small\Math{\PermRep{22275}\cong\IrrRep{1}\oplus\IrrRep{22}\oplus\underbrace{\IrrRep{252}\oplus\IrrRep{252}}_{}\oplus\underbrace{\IrrRep{1750}\oplus\IrrRep{1750}}_{}
\oplus\IrrRep{3520}\oplus\IrrRep{5103}\oplus\IrrRep{9625}}}\\[-6pt]
Time \textbf{C}: \Math{23} sec. Time \textbf{\textbf{Maple}}: \Math{11} sec.
	\item 
~\\[-17pt]
\textbf{66825}-dimensional representation of \Math{\mathbf{3.McL}}\\
Rank: \Math{14}. Suborbit lengths: \Math{1^3, 630, 2240^3, 5040^3, 8064^3, 20160}.\\[1pt]
{\small\Math{\PermRep{66825}\cong\IrrRep{1}\oplus\IrrRep{252}\oplus\IrrRep{252}\oplus\IrrRep{1750}\oplus\IrrRep{2772}\oplus\overline{\IrrRep{2772}}
\oplus\IrrRep{5103_\beta}\oplus\overline{\IrrRep{5103_\beta}}
\oplus\IrrRep{5103_\alpha}}\\
\Math{\hspace{48pt}\oplus\IrrRep{5544}\oplus\IrrRep{6336}\oplus\overline{\IrrRep{6336}}
\oplus\IrrRep{8064}\oplus\overline{\IrrRep{8064}}\oplus\IrrRep{9625}}}\\[1pt]
Time \textbf{C}: \Math{8} min \Math{45} sec. Time \textbf{\textbf{Maple}}: \Math{12} min \Math{59} sec.
	\item 
~\\[-17pt]
\textbf{22275}-dimensional representation \textbf{(a)} of \Math{\mathbf{McL\rtimes2}}\\
Rank: \Math{11}. Suborbit lengths: \Math{1, 112, 210, 420, 1120, 1260, 2520^2, 3360, 4032, 6720}.\\[1pt]
{\small\Math{\PermRep{22275}\cong\IrrRep{1}\oplus\IrrRep{22}\oplus\underbrace{\IrrRep{252}\oplus\IrrRep{252}}_{}\oplus\IrrRep{1750_\alpha}\oplus\IrrRep{1750_\beta}\oplus\IrrRep{3520}
\oplus\IrrRep{5103}\oplus\IrrRep{9625}}}\\[-5pt]
Time \textbf{C}: \Math{23} sec. Time \textbf{\textbf{Maple}}: \Math{5} sec.
\end{enumerate}
~\\[-45pt]
\paragraph{\textbf{Suzuki group} \Math{\mathbf{Suz}}.}\label{Suz}
\Math{\Ord=448345497600=2^{13}\cdot3^7\cdot5^2\cdot7\cdot11\cdot13,~ \mathrm{M}\cong\CyclG{6},~ \mathrm{Out}\cong\CyclG{2}.}
~\\[-19pt]
\begin{enumerate}
	\item 
\textbf{32760}-dimensional representation of \Math{\mathbf{Suz}}\\
Rank: \Math{6}. Suborbit lengths: \Math{1, 891, 1980, 2816, 6336, 20736}.\\[1pt]
{\small\Math{\PermRep{32760}\cong\IrrRep{1}\oplus\IrrRep{143}\oplus\IrrRep{364}\oplus\IrrRep{5940}\oplus\IrrRep{12012}\oplus\IrrRep{14300}}}\\[1pt]
Time \textbf{C}: \Math{54} sec. Time \textbf{\textbf{Maple}}: \Math{2} sec.
	\item 
~\\[-17pt]	
\textbf{65520}-dimensional representation of \Math{\mathbf{2.Suz}}\\
Rank: \Math{10}. Suborbit lengths: \Math{1^2, 891^2, 2816^2, 3960, 12672, 20736^2}.\\[1pt]
{\small\Math{\PermRep{65520}\cong\IrrRep{1}\oplus\IrrRep{143}\oplus\IrrRep{364_\alpha}\oplus\IrrRep{364_\beta}\oplus\overline{\IrrRep{364_\beta}}\oplus\IrrRep{5940}\oplus\IrrRep{12012}\oplus\IrrRep{14300}\oplus\IrrRep{16016}\oplus\overline{\IrrRep{16016}}}}\\[1pt]
Time \textbf{C}: \Math{6} min \Math{9} sec. Time \textbf{Maple}: \Math{11} sec.
	\item
~\\[-17pt]	
\textbf{98280}-dimensional representation of \Math{\mathbf{3.Suz}}\\
Rank: \Math{14}. Suborbit lengths: \Math{1^3, 891^3, 2816^3, 5940, 19008, 20736^3}.\\[1pt]
{\small\Math{\PermRep{98280}\cong\IrrRep{1}\oplus\IrrRep{78}\oplus\overline{\IrrRep{78}}\oplus\IrrRep{143}\oplus\IrrRep{364}\oplus\IrrRep{1365}\oplus\overline{\IrrRep{1365}}\oplus\IrrRep{4290}\oplus\overline{\IrrRep{4290}}\oplus\IrrRep{5940}\oplus\IrrRep{12012}}\\
\Math{\hspace{48pt}\oplus\IrrRep{14300}\oplus\IrrRep{27027}\oplus\overline{\IrrRep{27027}}}}\\[1pt]
Time \textbf{C}: \Math{57} min \Math{58} sec. Time \textbf{\textbf{Maple}}: \Math{6} min \Math{42} sec.\\[1pt]
\emph{Remark.} The \textit{\textbf{PreparePolynomialData}} program uses more than \Math{17.6} GB of memory for this task, which goes beyond the RAM of our PC, slowing down the calculations.
	\item 
\textbf{1782}-dimensional representation of \Math{\mathbf{Suz\rtimes{2}}}\\
Rank: \Math{3}. Suborbit lengths: \Math{1, 416, 1365}.\\[1pt]
{\small\Math{\PermRep{1782}\cong\IrrRep{1}\oplus\IrrRep{780}\oplus\IrrRep{1001}}}\\[1pt]
Time \textbf{C}: \Math{<1} sec. Time \textbf{\textbf{Maple}}: \Math{<1} sec.
	\item 
\textbf{5346}-dimensional representation of \Math{\mathbf{3.Suz\rtimes{2}}}\\
Rank: \Math{5}. Suborbit lengths: \Math{1, 2, 416, 832, 4095}.\\[1pt]
{\small\Math{\PermRep{5346}\cong\IrrRep{1}\oplus\IrrRep{132}\oplus\IrrRep{780}\oplus\IrrRep{1001}\oplus\IrrRep{3432}}}\\[1pt]
Time \textbf{C}: \Math{1} sec. Time \textbf{\textbf{Maple}}: \Math{<1} sec.
\end{enumerate}
~\\[-40pt]
\subsubsection{\textbf{Monster sections.}}\label{Monster} 
The main properties of the \textbf{Held group} \Math{\mathbf{He}} and the results of calculations for its representation of dimension \textbf{8330} are given in Section \ref{example}.
\begin{enumerate}
	\item 
\textbf{29155}-dimensional representation of \Math{\mathbf{He}}\\
Rank: \Math{12}. Suborbit lengths: \Math{1, 90, 120, 384, 960^2, 1440, 2160, 2880^2, 5760, 11520}.\\[1pt]
{\small\Math{\PermRep{29155}\cong\IrrRep{1}\oplus\IrrRep{51}\oplus\overline{\IrrRep{51}}\oplus\IrrRep{680}
\oplus\underbrace{\IrrRep{1275}\oplus\IrrRep{1275}}_{}\oplus\IrrRep{1920}\oplus\IrrRep{4352}\oplus\IrrRep{7650}\oplus\IrrRep{11900}}}\\[-5pt]
Time \textbf{C}: \Math{42} sec. Time \textbf{\textbf{Maple}}: \Math{11} sec.
	\item 
\textbf{8330}-dimensional representation of \Math{\mathbf{He\rtimes2}}\\
Rank: \Math{6}. Suborbit lengths: \Math{1, 105, 720, 1344, 1680, 4480}.\\[1pt]
{\small\Math{\PermRep{8330}\cong\IrrRep{1}\oplus\IrrRep{102}\oplus\IrrRep{680}\oplus\IrrRep{1275}\oplus\IrrRep{1920}\oplus\IrrRep{4352}}}\\[1pt]
Time \textbf{C}: \Math{3} sec. Time \textbf{\textbf{Maple}}: \Math{1} sec.
\end{enumerate}
~\\[-40pt]
\paragraph{\textbf{Fischer group} \Math{\mathbf{Fi_{22}}}.}\label{Fi22}
\Math{\Ord=64561751654400=2^{17}\cdot3^9\cdot5^2\cdot7\cdot11\cdot13,~ \mathrm{M}\cong\CyclG{6},~ \mathrm{Out}\cong\CyclG{2}.}
\begin{enumerate}
	\item 
\textbf{61776}-dimensional representation of \Math{\mathbf{Fi_{22}}}\\
Rank: \Math{4}. Suborbit lengths: \Math{1, 1575, 22400, 37800}.\\[1pt]
{\small\Math{\PermRep{61776}\cong\IrrRep{1}\oplus\IrrRep{3080}\oplus\IrrRep{13650}\oplus\IrrRep{45045}}}\\[1pt]
Time \textbf{C}: \Math{10} min \Math{6} sec. Time \textbf{\textbf{Maple}}: \Math{3} sec.
	\item 
\textbf{28160}-dimensional representation of \Math{\mathbf{2.Fi_{22}}}\\
Rank: \Math{5}. Suborbit lengths: \Math{1^2, 3159^2, 21840}.\\[1pt]
{\small\Math{\PermRep{28160}\cong\IrrRep{1}\oplus\IrrRep{352}\oplus\IrrRep{429}\oplus\IrrRep{13650}\oplus\IrrRep{13728}}}\\[1pt]
Time \textbf{C}: \Math{39} sec. Time \textbf{\textbf{Maple}}: \Math{2} sec.
	\item 
\textbf{56320}-dimensional representation of \Math{\mathbf{2.Fi_{22}\rtimes2}}\\
Rank: \Math{9}. Suborbit lengths: \Math{1^2, 728, 1080^2, 3159^2, 21840, 25272}.\\[1pt]
{\small\Math{\PermRep{56320}\cong\IrrRep{1}\oplus\IrrRep{1'}\oplus\IrrRep{352}\oplus\overline{\IrrRep{352}}\oplus\IrrRep{429_+}\oplus\IrrRep{429_-}\oplus\IrrRep{13650_+}\oplus\IrrRep{13650_-}\oplus\IrrRep{27456}}}\\[1pt]
Time \textbf{C}: \Math{3} min \Math{20} sec. Time \textbf{\textbf{Maple}}: \Math{5} sec.
\end{enumerate}
~\\[-40pt]
\paragraph{\textbf{Fischer group} \Math{\mathbf{Fi_{23}}}.\hspace*{-9pt}}\label{Fi23}
\Math{\Ord=4089470473293004800=2^{18}\cdot3^{13}\cdot5^2\cdot7\cdot11\cdot13\cdot17\cdot23,~\mathrm{M}\cong1,~ \mathrm{Out}\cong1.}\\
\textbf{31671}-dimensional representation of \Math{\mathbf{Fi_{23}}}\\
Rank: \Math{3}. Suborbit lengths: \Math{1, 3510, 28160}.\\[1pt]
{\small\Math{\PermRep{31671}\cong\IrrRep{1}\oplus\IrrRep{782}\oplus\IrrRep{30888}}}\\[1pt]
Time \textbf{C}: \Math{52} sec. Time \textbf{\textbf{Maple}}: \Math{1} sec.
~\\[-25pt]
\subsubsection{\textbf{Pariahs.}}\label{Pariahs} 
~\\[-30pt]
\paragraph{\textbf{Janko group} \Math{\mathbf{J_1}}.}
 \Math{\Ord=175560=2^{3}\cdot3\cdot5\cdot7\cdot11\cdot19,~ \mathrm{M}\cong1,~ \mathrm{Out}\cong1.}\\
\textbf{1045}-dimensional representation of \Math{\mathbf{J_1}}\\
Rank: \Math{11}. Suborbit lengths: \Math{1, 8, 28, 56^3, 168^5}.\\[1pt]
{\small\Math{\PermRep{1045}\cong\IrrRep{1}\oplus\IrrRep{56_+}\oplus\IrrRep{56_-}\oplus\IrrRep{76}\oplus\IrrRep{77_+}\oplus\IrrRep{77_-}      
\oplus\IrrRep{120_\alpha}\oplus\IrrRep{120_\beta} \oplus\IrrRep{120_\gamma}\oplus\IrrRep{133}\oplus\IrrRep{209}}}\\[1pt]
Time \textbf{C}: \Math{<1} sec. Time \textbf{\textbf{Maple}}: \Math{22} sec.
~\\[-25pt]
\paragraph{\textbf{Janko group} \Math{\mathbf{J_3}}.}
\Math{\Ord\farg{J_3}=50232960=2^{7}\cdot3^5\cdot5\cdot17\cdot19,~ \mathrm{M}\farg{J_3}\cong\CyclG{3},~ \mathrm{Out}\farg{J_3}\cong\CyclG{2}.}
\begin{enumerate}
	\item 
\textbf{14688}-dimensional representations \textbf{(a)} and \textbf{(b)} of \Math{\mathbf{J_3}}\\
Rank: \Math{14}. Suborbit lengths: \Math{1, 285, 342, 380, 570^2, 855^2, 1140^2, 1710^3, 3420}.\\[1pt]
{\small\Math{\PermRep{14688}\cong\IrrRep{1}\oplus\IrrRep{85}\oplus\overline{\IrrRep{85}}\oplus\underbrace{\IrrRep{1140}\oplus\IrrRep{1140}}_{}
\oplus\IrrRep{1215_+}\oplus\IrrRep{1215_-}\oplus\IrrRep{1615}\oplus\IrrRep{1920_\alpha}\oplus\IrrRep{1920_\beta}\oplus\IrrRep{1920_\gamma}\oplus\IrrRep{2432}}}\\[-5pt]
Time \textbf{C}: \Math{11} sec. Time \textbf{\textbf{Maple}}: \Math{1} min \Math{52} sec.\\
\emph{Remark.} 
\textsc{Atlas} \cite{atlas} contains two non-equivalent \Math{14688}-dimensional representations of \Math{J_3}, (a) and (b), which have the same decomposition structure. The differences are manifested in explicit expressions for irreducible projectors (and in the structure of orbitals).
	\item 
\textbf{6156}-dimensional representation of \Math{\mathbf{J_3\rtimes2}}\\
Rank: \Math{7}. Suborbit lengths: \Math{1, 85, 120, 510, 680, 2040, 2720}.\\[1pt]
{\small\Math{\PermRep{6156}\cong\IrrRep{1}\oplus\IrrRep{324}\oplus\IrrRep{646}\oplus\IrrRep{1140}\oplus\IrrRep{1215_+}\oplus\IrrRep{1215_-}\oplus\IrrRep{1615}}}\\[1pt]
Time \textbf{C}: \Math{1} sec. Time \textbf{\textbf{Maple}}: \Math{1} sec.
\end{enumerate}
~\\[-40pt]
\paragraph{\textbf{Rudvalis group} \Math{\mathbf{Ru}}.}
\Math{\Ord=145926144000=2^{14}\cdot3^3\cdot5^3\cdot7\cdot13\cdot29,~ \mathrm{M}\cong\CyclG{2},~ \mathrm{Out}\cong1.}
\begin{enumerate}
	\item 
\textbf{4060}-dimensional representation of \Math{\mathbf{Ru}}\\
Rank: \Math{3}. Suborbit lengths: \Math{1, 1755, 2304}.\\[1pt]
{\small\Math{\PermRep{4060}\cong\IrrRep{1}\oplus\IrrRep{783}\oplus\IrrRep{3276}}}\\[1pt]
Time \textbf{C}: \Math{<1} sec. Time \textbf{\textbf{Maple}}: \Math{<1} sec.
	\item 
\textbf{16240}-dimensional representation of \Math{\mathbf{2.Ru}}\\
Rank: \Math{9}. Suborbit lengths: \Math{1^4, 2304^4, 7020}.\\[1pt]
{\small\Math{\PermRep{16240}\cong\IrrRep{1}\oplus\IrrRep{28}\oplus\overline{\IrrRep{28}}\oplus\IrrRep{406}\oplus\IrrRep{783}\oplus\IrrRep{3276}
\oplus\IrrRep{3654}\oplus\IrrRep{4032}\oplus\overline{\IrrRep{4032}}}}\\[1pt]
Time \textbf{C}: \Math{12} sec. Time \textbf{\textbf{Maple}}: \Math{2} sec.
\end{enumerate}
~\\[-40pt]
\section{Concluding remarks}\label{concl}
For \textbf{\textit{PreparePolynomialData}}, the main limiting parameter is the representation dimension.
Our PC with 16 GB of RAM copes with dimensions not exceeding 100,000.
We can expect that with enough RAM, the program will cope with dimensions up to several hundred thousand.
\par
The main bottleneck of \textbf{\textit{SplitRepresentation}} is that it is based on the polynomial algebra methods, which are intrinsically algorithmically difficult.
The number of polynomial variables is equal to the rank \Math{\baseformN} of the representation to be split.
In practice, the program confidently splits representations with \Math{\baseformN\leq17}, although there are some examples with ranks 18 and 19.
However, representations of finite groups often have low ranks.
In particular, in \textsc{Atlas} \cite{atlas}, \Math{761} out of \Math{886}, or \Math{86\%}, permutation representations satisfy the condition \Math{\baseformN\leq17}.

~\\[-35pt]
\ack
I thank Yu.A. Blinkov, V.P. Gerdt, N.N. Vassiliev and R.A. Wilson for helpful discussions.
\end{document}